\documentclass[12 pt]{article} 
\usepackage[english]{babel}
\usepackage[active]{srcltx}
\usepackage{amsmath}
\usepackage{amsfonts,amssymb}
\usepackage[mathcal]{eucal}

\newcommand{\epsi}{\varepsilon}

\newcommand{\upartial}{\partial}
\newcommand{\rav}{\stackrel{\triangle}{=}}
\newcommand{\ravref}[1]{\stackrel{(\ref{#1})}{=}}
\newcommand{\leqref}[1]{\stackrel{(\ref{#1})}{\leq}}

\newcommand{\rref}[1]{$(\ref{#1})$}
\newcommand{\mm}[1]{{\mathbb{#1}}}
\newcommand{\ct}[1]{{\mathcal{#1}}}

\newcommand{\td}[1]{\widetilde{#1}}
\newcommand{\doc}{{\em{Proof.}}}
\newcommand{\bo}{\hfill {$\Box$}}
\newtheorem{definition}{Definition}
\newtheorem{remark}{Remark}
\newtheorem{theorem}{Theorem}
\newtheorem{lemma}{Lemma}
\newtheorem{corollary}{Corollary}

\begin{document}

\title{ON HAMILTONIAN AS LIMITING GRADIENT\\ IN  INFINITE HORIZON PROBLEM\thanks{
               Krasovskii Institute,, Yekaterinburg, Russia;\ \  Ural Federal University, Yekaterinburg, Russia 
}
}

\author{Dmitry Khlopin\\
{\it khlopin@imm.uran.ru} 
}

\maketitle

\begin{abstract}
Necessary conditions of optimality in the form of the Pontryagin Maximum Principle are derived for
the Bolza-type discounted problem with  free right end.
  The optimality is understood in the sense of the uniformly overtaking optimality. Such process is assumed to exist, and the corresponding   payoff of the optimal process (expressed in the form of improper integral) is assumed to converge in the Riemann sense.
    No other assumptions on the asymptotic behaviour of
    trajectories or adjoint variables
    are required.

In this paper, we prove that there exists a corresponding limiting solution of the Pontryagin Maximum Principle that satisfies the Michel transversality condition; in particular, the stationarity condition of the maximized Hamiltonian and the fact that  the maximized  Hamiltonian vanishes
at infinity are proved.  The connection of this condition with the limiting subdifferentials of payoff function
along the optimal process at infinity is showed.
The case of payoff without discount multiplier is also considered.


\end{abstract}

  \section*{Introduction}

      The main means of construction of necessary conditions of optimality for control problems is the Pontryagin Maximum Principle \cite{ppp}.
      In case of infinite horizon, the maximum principle is generally incomplete (see \cite{Halkin}): its relations offer no boundary condition at infinity. In absence of such  transversality
conditions, the PMP provides too many purportedly optimal solutions.

  Presently, many varieties of such conditions are constructed; a reference to all of them is not our intention. Nevertheless, let us note
\cite{kr_as,aucl,sch,Halkin,kam,Michel,pereira,smirn,ssbook}.
      One of such transversality conditions was proposed by Michel \cite{Michel}.
      If the dynamics of the equation is autonomous and the payoff is of the
form $\int_{0}^\infty e^{-rt}f_0(x,u)\,dt$, the condition may be rendered as
 $$-H[T]=
 \lambda^* r\int_{T}^{\infty} e^{-rt}f_0(x^*(t),u^*(t))\,dt\quad \forall T\geq0,$$
 where $(\psi^*,\lambda^*,x^*,u^*)$ satisfies the Pontryagin maximum principle.
Like the other transversality conditions for infinite horizon, it only becomes a necessary condition under additional assumptions, see
\cite[Sect. 6]{kr_as}.
There are many papers that prove the necessity of such  conditions under various assumptions. In \cite{ye}, the necessity was proved without assuming the dynamics to be smooth; in \cite{kam}, it was studied in the  calculus of variations setting; see \cite{olsil} for infinite horizon control problem with state constraints;
  in  \cite{ss} it was proved for the general statement, including the problems with fixed right end; in
\cite{kr_as}, under sufficiently weak assumptions on the summability, the connection of this condition with the Aseev--Kryazhimskii formula was studied. The assumptions used in this paper could not be embedded into assumptions of the above-mentioned papers; in particular, in contrast with \cite{kr_as,kam,olsil,ss}, here, as well as in
\cite{Michel}, the case of $\lambda^*=0$ is not generally excluded.

  Note that the Michel condition, if convenient, is only one-dimensional and, therefore,
 this condition, together with the core conditions of the maximum principle, can
 determine a unique solution candidate only for the problems with one state variable.
 In view of that, it is important to know not only when this  condition is necessary but also when it is consistent with other transversality conditions. For a similar analysis of the Aseev--Kryazhimskii formula, refer to \cite{kr_as}. Here, the Michel condition is used along with some
 limiting solution of the Pontryagin maximum principle (see \cite{Optim}); the limiting solution may be considered without assumptions on the asymptotic behaviour of
    trajectories or adjoint variables. The idea of the limiting solution can be traced to paper \cite{norv}; see its connection with the Aseev--Kryazhimskii formula in \cite{JDCS}. The general case of Bolza-type infinite horizon problem with  free right end was studied in \cite{Optim}. In this paper, we prove the existence of a limiting solution of PMP that satisfies the Michel condition for uniformly overtaking optimal control; the arising transversality conditions are expressed in the form of  limiting gradients of payoff function at infinity.
 The proof itself combines the ideas from \cite{Michel} with the proof of the Pontryagin Maximum Principle from \cite{Optim}.

 The paper is structured as follows. First we describe the problem statement, impose the general conditions and propositions; at the same section, we provide the required definitions from the smooth analysis.
 In Section~\ref{187}, in addition to the PMP relations and definition of a limiting solution to the Pontryagin maximum principle, we specify the computation of limiting gradients of the payoff function at infinity. In the next section, we formulate the main result (Theorem~\ref{4}) and a number of its simple corollaries.
 The last two items contain, respectively, the preliminary lemmas and the proof of Theorem~\ref{4}.

  \section{Problem statement and definitions}

  We consider the time interval
  ${\mathbb{T}}\rav
  {\mathbb{R}}_{\geq 0}.$
The phase space of the control system is the finite-dimensional Euclidean space ${\mathbb{X}}\rav{\mathbb{R}}^m$.

  Consider the following optimal control problem
\begin{subequations}
\begin{eqnarray}
    \textrm{Minimize } l(b)+\int_{0}^\infty e^{-rt} f_0(x,u) dt  \label{sys0_}\\
    \textrm{subject to } \dot{x}=f(x,u),\quad u\in U,  \label{sys_}\\
    x(0)\in\ct{C}.    \label{sysK_}
\end{eqnarray}
\end{subequations}
  Here, the function $f_0$ is scalar; $x$ is the  state variable taking values in  ${\mathbb{X}};$ and $u$ is the control parameter.

Suppose that  $U$ is
a Borel subset of a finite-dimensional Euclidean space. As for the class of admissible controls, we consider the set of measurable functions $u(\cdot)$ bounded for any time compact
 such that $u(t)\in U$ holds for a.a. $t\in{\mathbb{T}}$.
 Denote the set of admissible controls
   by~${\ct{U}}$.

 We assume the following conditions hold:
\begin{itemize}
  \item  ${\mathcal{C}}$ is a closed subset of ${\mathbb{X}}$;
  \item   $l$ is taken to be locally Lipshitz continuous on $x$;
  \item  $f$ is Borel measurable in $u$ and continuously differentiable
in $x$;
  \item   for each admissible control $u$, the map  $(t,x)\mapsto f(x,u(t))$ satisfies the
 sublinear growth condition (see, for example, \cite[1.4.4]{tovst});
  \item   $f_0$ is measurable in $u$, continuously differentiable
in $x$,  and lower semicontinuous in $u$;
  \item  $\frac{\partial f}{\partial x},\frac{\partial f_0}{\partial x}$ are measurable in $u$ and
 locally Lipshitz continuous on  $x$.
\end{itemize}

For each admissible control~$u$, 
and position $b\in{\mathbb{X}}$,
we can consider a solution of \rref{sys_} for $x(0)=b.$
The solution is unique and it can be extended to the whole ${\mathbb{T}}$.
Let
us denote it by~$x(b,u;\cdot).$

The pair $(x,u)$ will be called an admissible control process
if~$u\in {\mathfrak{U}},$ $x(0)\in{\mathcal{C}},$ $x(\cdot)=x(x(0),u;\cdot).$

\begin{definition}
 If an admissible process $(x^*,u^*)$ satisfies
\begin{eqnarray*}
 \limsup_{T\to\infty}
  \bigg[
    l(x^*(0))&+&\int_0^T e^{-rt}f_0\big(x^*(t),u^*(t)\big)\,dt
    -\\
      \inf_{(b,u)\in{\mathcal{C}}\times{\mathcal{U}}}\Big(l(b)&+&\int_0^T e^{-rt}f_0\big(x(b,u;t),u(t)\big)dt\Big)
     \bigg]
      \leq 0,
\end{eqnarray*}
call it
    a uniformly overtaking optimal process for \rref{sys0_}--\rref{sysK_}.
\end{definition}

  Hereinafter assume there exists an optimal uniformly overtaking process $(x^*,u^*)$.
  Set $b_*=x^*(0).$ We are not going to impose any conditions that guarantee the existence of such a solution; for various existence theorems, refer to, for example,
  \cite{baum,baumeister,slovak,z1}.

      Let the improper integral
  $$\int_{0}^{\infty} e^{-rt} f_0(x^*(t),u^*(t))\,dt$$
  converge in the Riemann sense, i.e.,
  \begin{equation}
  \label{raz_}
     J_{**}\rav\lim_{t\to\infty}\int_{0}^{t} e^{-rt}f_0(x^*(t),u^*(t))\,dt\in\mm{R}.
  \end{equation}

Let us now define scalar functions $J^{0},\bar{J}^0$ by the following rule: for all $b\in\mm{X},T,s\geq 0$,
\begin{eqnarray*}
   J^0(b,s;T)&\rav&\int_{0}^{T} e^{-r(t+s)}f_0\big(x(b,u^*;t),u(t)\big)\,dt.\\
   \bar{J}^{0}(b;T)\rav J^0(b,0;T)&=&\int_{0}^{T} e^{-rt}f_0\big(x(b,u^*;t),u(t)\big)\,dt.
\end{eqnarray*}
To continue, we need to define subgradients of these payoffs  at infinity.
To this end, let us introduce the necessary notions of  convex analysis \cite{cl_new},\cite[Sect.4]{vinter}.

Consider a finite-dimensional Euclidian space  $E$, and a lower semicontinuous function
$g:E\to{\mathbb{R}}\cup\{+\infty\}$.
 A vector $\zeta\in E$
is said to be a proximal subgradient of $g$  at $b\in E$ if
there exist a
neighborhood $\Omega$ of $b$ and a number $\sigma\geq 0$ such that
$g(\xi) \geq g(b) + \zeta(\xi-b)-\sigma||\xi-b||^2$
 for all
 $\xi\in \Omega.$
The set of proximal subgradients at $b$  is denoted $\upartial_P g(b)$,
and is referred to as the proximal subdifferential. This set  is nonempty for all $b$ in
a dense subset of $\{b\,|\,g(b)<+\infty\}.$
 Following
 \cite[Theorem 4.6.2(a)]{vinter},
  denote the limiting subdifferential of $g$ at $b$ by $\upartial_L g(b)$; it
 consists of all
$\zeta$ in $E$  such that
\[
\exists \textrm{ sequences of }  y_n\in{\mathbb{X}},\zeta_n\in \upartial_P g(y_n),y_n\to b,
\zeta_n\to\zeta.\]
 Following
 \cite[Theorem 4.6.2(b)]{vinter},
 denote the singular limiting (asymptotic limiting) subdifferential of $g$ at $b$ by $\upartial^0_{L} g(b)$; it consists of all
$\zeta$ in $E^*$  such that
\[
\exists \textrm{ sequences of }  y_n\in E,\lambda_n\in{\mathbb{T}},\zeta_n\in \upartial_P g(y_n),y_n\to b,\lambda_n\downarrow 0,
\lambda_n\zeta_n\to\zeta.\]
If $g$ is Lipshitz continuous near $b$, then $\upartial_{L} g(b)$ is nonempty, moreover $co\, \upartial_L g(b)=\upartial_{Clarke} g(b),\upartial^0_{L} g(b)=\{0\}$
(see \cite[Sect. 4]{vinter}).

Following the same idea, define the subgradients of $g$ at infinity, or, more accurately, along on arbitrary unboundedly increasing sequence  of positive  $\tau.$ Fix a sequence $\tau.$

Denote ${\mathcal{T}}\rav\{\tau_n\,|\,n\in{\mathbb{N}}\}.$
For a differentiable function $g:E\times{\mathbb{T}}\to{\mathbb{X}}$,
similarly to the definitions of limiting subdifferential and singular limiting subdifferential,
let us introduce the generalized subdifferential of $g$ at the infinite point $(b,\infty_\tau)$, or rather at $b$ with infinity along $\tau$, by the following rule:
\begin{eqnarray*}
\upartial^{1}_L g(b,\infty_\tau)&\rav&
  \{\zeta\,|\,\exists \textrm{ sequences of }  y_n\in E,t_n\in{\mathcal{T}},
  \zeta_n\in \upartial_P g(y_n,t_n),\\ &\ &y_n\to b,
 t_n\to\infty,\zeta_n\to\zeta\}.
\end{eqnarray*}
  Since in the general case it may be empty, let us also introduce a
  singular subdifferential in the following way:
\begin{eqnarray*}\upartial^{0}_L g(b,\infty_\tau)&\rav&
  \{\zeta\,|\,\exists \textrm{ sequences of }  y_n\in E,t_n\in{\mathcal{T}},\lambda_n\in{\mathbb{T}},\zeta_n\in \upartial_P g(y_n,t_n),\\
  &\ &y_n\to b,t_n\to\infty,\lambda_n\to 0,
\lambda_n\zeta_n\to\zeta\}.
\end{eqnarray*}

 Note that the mentioned definitions can be rewritten otherwise.
 First of all, in the last two definitions, $\upartial_P g(b,t_n)$ can be replaced with $\upartial_L g(b,t_n)$ because every element from
 $\upartial_L g(b,t_n)$ can be approximated with arbitrary precision by an element from $\upartial_L g(y,t_n)$ for some $y$ that is arbitrarily close to $b$. Moreover, in the definition of $\upartial^{1}_L g$, one can replace
 $\zeta_n\to\zeta$ with $\lambda_n\zeta_n\to\zeta$ under the condition $\lambda_n\to 1.$
 Thus we obtain the equivalent form:
\begin{eqnarray*}
  \upartial^{\lambda}_L g(b,\infty_\tau)&=&
  \{\zeta\,|\,\exists \textrm{ sequences of }  y_n\in E,t_n\in{\mathcal{T}},\lambda_n\in{\mathbb{T}},\zeta_n\in \upartial_L^1 g(y_n,t_n),\\
  &\ &y_n\to b,t_n\to\infty,\lambda_n\to \lambda,
\lambda_n\zeta_n\to\zeta\}\quad \forall\lambda\in\{0,1\}.
\end{eqnarray*}

Remember that the limiting normal cone $N_L^{{\mathcal{C}}}(b)$ of ${\mathcal{C}}$ at $b$
is the  limiting subdifferential $\upartial^{1}_L\delta_{\mathcal{C}}(b)$ of the indicator function $\delta_{\mathcal{C}}$
of the set ${\mathcal{C}}$ (see, for example, \cite[Proposition 1.18]{kruger}).

\section{Limiting solution of the Pontryagin Maximum Principle}
\label{187}

   Let us now proceed to the relations of the Pontryagin Maximum Principle.

  Let the Hamilton--Pontryagin function
  $H:{\mathbb{X}}\times {U}\times{\mathbb{X}}\times{\mathbb{T}}\times {\mathbb{T}}\mapsto{\mathbb{R}}$
  be given by
   $H(x,u,\psi,\lambda,t)\rav\psi f\big(x,u\big)-\lambda
   e^{-rt}f_0\big(x,u\big).$
 Let us introduce the relations:
\begin{subequations}
 \begin{eqnarray}
   \label{sys_x}
       \dot{x}(t)&=& f\big(x(t),u(t)\big);\\
   \label{sys_psi}
       -\dot{\psi}(t)&=&\frac{\partial
       H}{\partial x}\big(x(t),u(t),\psi(t),\lambda,t\big);\\
   \label{maxH}
\sup_{u'\in
U(t)}H\big(x(t),u',\psi(t),\lambda,t\big)&=&
        H\big(x(t),u(t),\psi(t),\lambda,t\big);
   \end{eqnarray}
 for norming, it would also be convenient to use one of the following conditions:
 \begin{eqnarray}
   \label{dob}
   ||\psi(0)||+\lambda=1;
        \\
   \label{dob01}
   \lambda\in\{0,1\}.
   \end{eqnarray}
\end{subequations}
 It is easily seen that, for each $u\in{\ct{U}}$ for each initial
  condition, system~\rref{sys_x}--\rref{sys_psi} has a local solution,
   and each solution of these relations can be extended to the
    whole~${\mathbb{T}}$.
\begin{remark}
   Since the right-hand side of~\rref{sys_psi_k}--\rref{dob_k} is homogeneous by $(\psi,\lambda)$,
     a nontrivial solution of \rref{sys_x}--\rref{maxH} with \rref{dob} exists iff there exists a nontrivial solution
     of \rref{sys_x}--\rref{maxH} with \rref{dob01}.
\end{remark}

  Although the PMP holds for a rather general
  infinite-horizon control problem \cite{Halkin}, its system of relations  \rref{sys_x}--\rref{dob} is generally incomplete and requires an additional boundary condition. Many such conditions, which hold under various supplementary assumptions (imposed, first of all, on the asymptotic behavior of the adjoint variable) were offered
  (see the reviews in \cite{kr_as,ye}). In the general case, the result below does not
  require such assumptions \cite{JDCS}:
\begin{theorem}
    Let the process $(x^*,u^*)$ be  a uniformly overtaking process for problem \rref{sys0_}--\rref{sysK_} with singleton $\ct{C}$. Let $\tau$ be an unbounded increasing sequence of positive numbers.

    Then,
    for  $(x^*,u^*)$ there exists a $\tau$-limiting solution $(\psi^*,\lambda^*)$ of  system \rref{sys_x}--\rref{maxH} satisfying \rref{dob}.
\end{theorem}



\begin{definition}
    A nontrivial solution
   $(\lambda^*,\psi^*)$ of  \rref{sys_x}--\rref{maxH}
   associated with $(x^*,u^*)$ is called $\tau$-limiting (or just limiting)
   if,   for some subsequence
                     $\tau'\subset\tau$,
                      $(x^*,\psi^*,\lambda^*)$ is a pointwise limit
      of
        solutions $(x_n,\psi_n,\lambda_n)$ of the boundary value problems
\begin{subequations}
 \begin{eqnarray}
   \label{sys_x_k}
       \dot{x}_n(t)&=& f\big(x_n(t),u^*(t)\big);\\
        \label{sys_psi_k}
       -\dot{\psi}_n(t)&=&\frac{\partial
       H}{\partial x}\big(x_n(t),u^*(t),\psi_n(t),\lambda_n,t\big);\\
   \label{sys_l_k}
       \dot{\lambda}_n(t)&=&0;\\
   \label{dob_k}
   \psi_n(\tau'_n)&=&0 
 \end{eqnarray}
 \end{subequations}
 on the interval $[0,\tau'_n].$
\end{definition}
\begin{remark}
    Without loss of generality we can say  that
      $\lambda_n+||\psi_n(0)||=\lambda^*+||\psi^*(0)||=1.$ Or, if $\lambda^*>0$, then we can say that
      $\lambda_n=\lambda^*=1$.
\end{remark}

This definition of $\tau$-limiting solution $(\psi^*,\lambda^*)$ of  system \rref{sys_x}-\rref{dob}
has several equivalent formulations.

First of all,
let us use the fact that system
 \rref{sys_psi_k}--\rref{sys_l_k}  is linear. Denote by
${\mathbb{L}}$ the linear space of all real $m\times m$ matrices; here $m=dim\,{\mathbb{X}}$.
For each $\xi\in{\mathbb{X}}$, there exists a solution ${A}(\xi;t)\in C({\mathbb{T}}, {\mathbb{L}})$ of the Cauchy problem
 \begin{equation*}
 \frac{d{A}(\xi;t)}{dt} =\frac{\partial f }{\partial x}
 \big(x(\xi,u^*;t),u^*(t)\big)
  A(\xi;t),\quad A(\xi;0)=1_{\mathbb{L}}.
\end{equation*}
Define the vector-valued function $I$ of time by the following rule: for every
  $T\in{\mathbb{T}}$,
  \begin{equation*}
  I(\xi;T)\rav\int_0^T
   e^{-rt}\frac{\partial f_0}{\partial x}
    \big(x(\xi,u^*;t),u^*(t)\big)
\, A(\xi;t)
  \,dt.
\end{equation*}

  Now, a solution $(x_n,\psi_n,\lambda_n)$ of system \rref{sys_x_k}--\rref{sys_l_k}
satisfies  the Cauchy formula:
  \begin{equation}
   \label{4A}
   \psi(t)=\big(\psi(0)+\lambda I(x(0);t)\big)A^{-1}(x(0);t)
   \qquad \forall t\in{\mathbb{T}}.
\end{equation}

  Note that thanks to $\psi_n(\tau_n)=0$,
  we have $\psi_n(t)=\lambda_n\big( I(x_n(0);t)-I(x_n(0);\tau_n)\big)A^{-1}(x_n(0);t);$
  in particular,
  \begin{equation}
   \label{4B}
  \psi_n(0)=-\lambda_n I(x_n(0);\tau_n).
\end{equation}

 Passing to the limit and using the expression for $I$ specified before, one can obtain the formulas for $\psi^*(0).$ In particular, if there exists a finite limit of
 $I(b;t)$ as $b\to b_*,t\to\infty,$ we see that the limiting co-state arc is unique up to a positive multiplication and the Aseev--Kryazhimskii formula holds:
  \begin{equation*}
  -\psi^*(0)=\int_0^\infty
   e^{-rt}\frac{\partial f_0}{\partial x}
    \big(x^*(t),u^*(t)\big)
\, A(b_*;t)
  \,dt,\qquad \lambda^*=1.
\end{equation*}
  Assumptions under which this expression is a necessary condition of optimality that are relatively easy to check may be found in \cite{kr_as,kab,JDCS}. This formula may not point towards a solution of the PMP even if the integral converges in the Lebesgue sense, see \cite{Optim}.  For details on the other (the more general formulas), see \cite{JDCS}.

  To make an all-encompassing formulas of limiting co-state arc one can use the terms of limiting subdifferentials of the payoff function $\bar{J}^0$ at infinity:
in \cite[Theorem~3.1]{Optim} it was proved that
 \begin{theorem}
{A solution} $(\psi^*,\lambda^*)$ {of} \rref{sys_x}--\rref{maxH},\rref{dob01} {associated with} $(x^*,u^*)$
      {is  $\tau$-limiting iff}
{a  solution} $(\psi^*,\lambda^*)$ {of} \rref{sys_x}--\rref{maxH},\rref{dob01} {associated with} $(x^*,u^*)$ is
nontrivial, and satisfies
$$\psi^*(0)\in \upartial^{\lambda^*}_{L} (-\bar{J}^{0})(x^*(0);u^*,\infty_\tau).$$
\end{theorem}

In following, we show that the condition \rref{raz_} implies the similar formula for the Hamiltonian, or rather the pair $(\psi,-H)$.

  For every $\vartheta>0$, define a control $u^{*-\vartheta}\in\ct{U}$ by the rule $u^{*-\vartheta}(t)=u^{*}(t+\vartheta).$
  Now, for every $b\in\mm{X}$, there exists $\xi\in\mm{X}$ such that $x(\xi,u^*;\vartheta)=b.$
  Then, $x(\xi,u^*;\vartheta+t)=x(b,u^{*-\vartheta};t)$ for all $t\geq 0.$
  We can now provide a definition valid for all
   $s\in\mm{R},T\geq \vartheta$:
\begin{eqnarray*}
   J^{\vartheta}(b,s;T)\rav
    J^{0}(\xi,s;T)-J^{0}(\xi,s;\vartheta)
    \quad
   \bar{J}^{\vartheta}(x_*;T)\rav {J}^{\vartheta}(x_*,0;T).
\end{eqnarray*}

  Note that, for all $T\in{\mathbb{T}}$,
  \begin{eqnarray}
   \frac{\partial J^0}{\partial b}(b,s;T)\equiv e^{-rs}I(b;T),\quad
   \frac{\partial J^0}{\partial s}(b,s;T)=-rJ^0(b,s;T)\label{414}.
  \end{eqnarray}
  Now, for every solution of system \rref{sys_x_k}--\rref{sys_l_k}, the following identities hold:
  \begin{eqnarray*}
   \frac{\partial J^\vartheta}{\partial s}(x_n(\vartheta),s;T)&=& -rJ^\vartheta(x_n(\vartheta),s;T);\\
   \frac{\partial J^\vartheta}{\partial b}(x_n(\vartheta),s;T)&=&
      \frac{\partial}{\partial b}\Big[
      J^0(x_n(0),s;T)-J^0(x_n(0),s;\vartheta)\Big]
      \bigg[\frac{\partial x(\xi,\vartheta;u^*)}{\partial \xi}\Big|_{\xi=x_n(0)}\bigg]^{-1}\\
      &\ravref{414}&
      e^{-rs}\big[I(x_n(0);T)-I(x_n(0);\vartheta)\big]A^{-1}(x_n(0);T)\\&\ravref{4A}&
      e^{-rs}\Big(\psi_n(T)A(x_n(0);T)A^{-1}(x_n(0);\vartheta)-\psi_n(\vartheta)\Big)/\lambda_n.
  \end{eqnarray*}
   Since all these mappings are continuous, for $T=\tau_n,s=0$, we obtain
\begin{eqnarray}
 \label{499}
     \upartial^1_L (-\bar{J}^t)(x_n(t),0;\tau_n)&=&\{\psi_n(t)/\lambda_n\};\\
 \label{500}
  \upartial^1_L (-J^t)(x_n(t),0;\tau_n)&=&\{(\psi_n(t)/\lambda_n,rJ^t(x_n(t),0;\tau_n))\}.
\end{eqnarray}
    Let us also note that $\lambda_n\to \lambda^*=0$ exactly when
     $-\psi_n(0)=\lambda_n I(x_n(0);\tau_n)\to-\psi^*(0),$  i.e., when $||I(x_n(0);\tau_n)||\to\infty.$


\section{The main result}

\begin{subequations}
\begin{theorem}
\label{4}
    Let the process $(x^*,u^*)$ be uniformly overtaking  optimal for problem \rref{sys0_}--\rref{sysK_}.
 Assume condition \rref{raz_} to hold.
Take an arbitrary unboundedly increasing sequence of times ~$\tau_n.$

    Then,
    for $(x^*,u^*)$ there exists a nontrivial $\tau$-limiting solution $(\psi^*,\lambda^*)$ of system \rref{sys_psi}-\rref{maxH} for $\lambda^*\in\{0,1\}$
    such that
    $$H(x^*(t),u^*(t),\psi^*(t),\lambda^*,t)=\psi^*(t)f(x^*(t),u^*(t))-\lambda^*e^{-rt}f_0(x^*(t),u^*(t)),$$
    for almost all $t$, coincides with  a continuous function $H^*:\mm{T}\to\mm{R}$,
    and,  for all $T\in\mm{T}$, the function $H^*$ satisfies
\begin{eqnarray}
\label{413} \textrm{(vanishing of Hamiltonian)} \qquad\quad\quad\ \lim_{t\to \infty} H^*[t]&=&0;\\
\label{412}\textrm{(stationarity condition)} \qquad\qquad\qquad\ \   -H^*[T]
 &=&\lambda^*r\lim_{n\to\infty}\bar{J}^T(x^*(T);\tau_n)\\
&=&\nonumber\lambda^* r\big(J_{**}-J^0(x^*(0),0;T)\big)\\
  &=&\nonumber\lambda^* r\int_{T}^{\infty} e^{-rt}f_0(x^*(t),u^*(t))\,dt;\\
\label{400}
\textrm{(transversality  condition at zero)} \qquad\qquad     \psi^*(0)&\in& \lambda^*\upartial_L l(x^*(0))+N_L^{{\mathcal{C}}}(x^*(0));\\
\textrm{(limiting   condition for x)} \qquad\qquad\qquad\quad
\label{409}
\psi^*(T)&\in& \upartial^{\lambda^*}_L (-\bar{J}^T)(x^*(T);\infty_{\tau});\\
\textrm{(limiting   condition for (x,t))} \quad \label{410}
(\psi^*(T),-H^*[T])&\in& \upartial^{\lambda^*}_L (-J^T)(x^*(T),0;\infty_{\tau}).
\end{eqnarray}
 Moreover, if $\lambda^*=0$, then $\psi^*(0)\neq 0$ holds, the sequence $I(x_n(0),\tau_n)$ is unbounded, and for almost all  $T>0$
\begin{eqnarray}
\label{403}
H^*[T]&\equiv&0\qquad\forall T\geq 0;\\
\label{404}
\psi^*(T)f(x^*(T),u^*(T))&=&0\qquad \forall \,a.a.\,T\geq 0.
\end{eqnarray}
\end{theorem}

Note that  if $f_0$ is bounded and $r>0$, then \rref{raz_} holds. Such assumption is used, for example, in \cite{ye}

Let us also make several simple observations.
 \begin{corollary}
Under assumptions of the theorem, let $r=0$; then, in addition to  \rref{413}--\rref{410}, we also have \rref{403} and
\begin{eqnarray}
\label{423}
\psi^*(T)f(x^*(T),u^*(T))=\lambda^*f_0(x^*(T),u^*(T))\qquad \forall \,a.a.\,T\geq 0.
\end{eqnarray}
 \end{corollary}

 Another one of them is about the converse of Hartwick’s rule in resource economics (see \cite{ss,WA})
 \begin{corollary}
Under assumptions of the theorem, let  $f_0(x^*(t),u^*(t))=C$ hold true for a certain constant $C$ for almost all  $t\geq 0.$
Then, there exists a nontrivial $\tau$-limiting solution $(\psi^*,\lambda^*)$ of system \rref{sys_psi}-\rref{maxH}
 for which, in addition to \rref{413}--\rref{410},  \rref{404} also holds.
 \end{corollary}
\end{subequations}
 {\doc}
Indeed, replace $f_0(x,u)$ with the function $f_C(x,u)\rav f_0(x,u)-C.$  The optimal process remains optimal, condition \rref{raz_} holds, and the solution of the PMP does not change. Apply the proved theorem. Then,
\rref{404} holds for almost all  $T>0$ by virtue of
\begin{eqnarray*}
H^*_C[T]\rav\psi^*(T)f(x^*(T),u^*(T))-\lambda^*e^{-rT}f_C(x^*(T),u^*(T))=\psi^*(T)f(x^*(T),u^*(T)),\\
 -H^*_C[T]\ravref{412}
 \lambda^* r\int_{T}^{\infty} e^{-rt}f_C(x^*(t),u^*(t))\,dt=0.
\end{eqnarray*}
\bo

 \begin{corollary}
\label{67}
 In some neighborhood
 $\Omega$ of the point $b_*$,  let the value function $V^T(b)$ of the problem
\begin{eqnarray*}
    \textrm{Minimize } \int_{0}^T f_0(x,u) dt  
    \\
    \textrm{subject to } \dot{x}=f(x,u),\quad u\in U  
    \\
    x(0)=b.    
\end{eqnarray*}
 be such that, for a Lipshitz function  $V^\infty$ defined in that neighborhood and a number
 $H^\infty\in\mm{R}$, we have
\begin{eqnarray}
     \lim_{T\to\infty} \big[V^T(b)+H^\infty T\big]=V^\infty(b) \quad \forall b\in\Omega,
     \label{dva_1}
\end{eqnarray}
 and this limit is uniform for $b\in\Omega.$

  In addition, let the control $u^*\in\ct{U}$ satisfy
\begin{eqnarray}
     \lim_{T\to\infty} \big[\bar{J}^0(b_*,u^*;T)+H^\infty T\big]=V^\infty(b_*).
     \label{raz_1}
\end{eqnarray}

 Then, for $(x^*,u^*)$, there exists a $\tau$-limiting solution  $(\psi^*,\lambda^*)$ of the PMP relations such that
 $\lambda^*=1$ and
\begin{eqnarray}
 \label{441}
 \psi^*(0)\in \upartial^1_L (-V)(b_*),&\ &\\
 \label{423_}
\psi^*(T)f(x^*(T),u^*(T))=f_0(x^*(T),u^*(T))+H^\infty&\ &\qquad \forall \,a.a.\,T\geq 0.
\end{eqnarray}
 \end{corollary}
 {\doc}
For every $T>0,b\in\Omega$, consider the problem
 \begin{eqnarray*}
    \textrm{Minimize } -V^T(b)+\int_{0}^T [f_0(x,u)+H^\infty] dt  
    \\
    \textrm{subject to } \dot{x}=f(x,u),\quad u\in U,\\  
    x(0)=b.    
\end{eqnarray*}
Without loss of generality we can assume  that $\Omega$ is closed.
 It is easy to see that the function of optimal value for this problem equals $H^\infty T$.
 Note that conditions \rref{raz_1} and \rref{dva_1} now imply condition \rref{raz_}, as well as the fact that $u^*$
 is uniformly overtaking optimal in the problem  \rref{sys0_}--\rref{sysK_} (with $r=0$,$\ct{C}=\Omega$)

  Then, the result of the theorem holds for it; in particular, there exists
    $\psi^*(0)\in \lambda^*\upartial_L^1 (-V)(b_*)+N_L^{\Omega}(b_*).$
    Since $N_L^{\Omega}(b_*)=\{0\}$, from $\lambda^*=0$ one would imply $\psi^*(0)=0$, which contradicts the nontriviality of the
    $\tau$-limiting solution. Then, $\lambda^*=1$ and \rref{441}.
 Writing out \rref{423} for this problem, we obtain \rref{423_}. \bo

 Note that condition \rref{441} is nothing else but the economical interpretation of a co-state arc as a shadow price.
 It is proved under varying assumptions on the system, for example, in \cite{a_new,aucl,sch,pereira,sagara}.

\section{Auxiliary lemmas.}

 Let $E,\Upsilon$ be a finite-dimensional Euclidean spaces.
  Consider a map   $a:{E}\times\Upsilon\times {\mathbb{T}}\mapsto E$.

As for the class of admissible controls, we consider any nonempty subset of measurable functions $\alpha(\cdot)$ bounded for any time compact
 such that $\alpha(t)\in \Upsilon$ holds for a.a. $t\in{\mathbb{T}}$.
 Denote the set of admissible controls
   by~${\ct{A}}$.

 For each admissible control
$\alpha\in{\ct{A}}$,  consider the differential equation:
\begin{equation}
   \label{a}
   \dot{y}=a(y(t),\alpha(t),t),\qquad
   \forall t\geq 0.
\end{equation}

   Assume that,
   for each admissible control $\alpha$,
 the map $(y,t)\mapsto a(y,\alpha(t),t)$ is a
 Carath\'eodory map;
 on each bounded subset, a map $(y,\alpha,t)\mapsto a(y,\alpha,t)$
 is integrally bounded and  locally Lipshitz
   continuous on $x$; moreover, each its local solution of \rref{a} can be extended onto the
whole~${\mathbb{T}}$ \cite{tovst}. For every admissible control $\alpha$, let us denote the
family of all solutions $y\in C({\mathbb{T}},E)$ of
system~\rref{a} by ${{\mm{Y}}}[\alpha]$.

 Consider any  admissible control $\alpha^*\in\ct{A}$ and a compact set ${\mathcal{S}}$.
Let us fix $\alpha^*,\mathcal{S}$.

 For every point
 $(y_*,\vartheta)\in{E}\times{\mathbb{T}}$, there exists a unique solution $y\in C({\mathbb{T}},E)$
 of the equation
\begin{equation}
   \label{1667}
    \dot{y}=a(y(t),\alpha^*(t),t),\quad y(\vartheta)=y_*.
\end{equation}
Let us denote its
initial position~$y(0)$ by $\varkappa(y_*,\vartheta)$.

In \cite{Optim} for such system~$a$
  with  the designated control $\alpha^*$ and the compact set ${\mathcal{S}}$,
 the map $$w:\Upsilon\times\Upsilon\times\mm{T}\to\mm{T}$$ was constructed. It has the following properties:
   \begin{itemize}
     \item    the mapping $\mm{T}\ni t\mapsto w(\alpha'(t),\alpha''(t),t)$ is Borel measurable for each $\alpha',\alpha''\in{\ct{A}}$;
     \item    the mapping $\Upsilon\times\Upsilon\ni(\alpha',\alpha'')\mapsto w(\alpha',\alpha'',t)$ is lower-semicontinuous for a.e. $t\in\mm{T}$;
     \item    for any $(\alpha',\alpha'',t)\in\Upsilon\times\Upsilon\times\mm{T}$, $w(\alpha',\alpha'',t)=0$ iff $\alpha'=\alpha''$;
     \item
   for any $(\alpha',\alpha'',t)\in\Upsilon\times\Upsilon\times\mm{T}$, $w(\alpha',\alpha'',t)=w(\alpha'',\alpha',t)\geq 1$ if $\alpha'\neq\alpha''$;
     \item
   the following lemmas hold (see \cite[Lemmae A.1--A.3]{Optim}):
   \end{itemize}
\begin{lemma}
\label{dop0}
 For every $T>0$, the mapping $${\ct{A}}\ni(\alpha',\alpha'')\mapsto\rho(\alpha',\alpha'',T)\rav \int_0^T w(\alpha'(t),\alpha''(t),t)dt$$ defines a metric on
\[{\ct{A}}_T\rav \{\alpha\in{\ct{A}}\,|\,\alpha(t)=\alpha^*(t) \ \forall t>T\};\]
 under this metric, the space ${\ct{A}}_T$ becomes a complete metric space, and the convergence in this metric is not weaker than the convergence in measure.

 In particular, if for some unbounded increasing sequence of times $\tau_n$, for some sequence of $\alpha_n\in{\ct{A}}_{\tau_n}$,
the sequence of
$\rho(\alpha^*,\alpha_n,\tau_n)$ tends to zero, then the sequence of $\alpha_n$ converges in the measure to $\alpha^*$ on
the whole ${\mathbb{T}}$.
 \end{lemma}
\begin{lemma}
\label{dop}
For arbitrary $\alpha\in{{\ct{A}}}$, $T>0$, every
 solution  $y\in{\mm{Y}}[\alpha], y(0)\in{\mathcal{S}}$
 of equation \rref{a} satisfies
\begin{equation}
\label{1000}
 \big|\big|\varkappa(y(t),t)-y(0)\big|\big|\leq \rho(\alpha^*,\alpha,t) \qquad \forall t\in [0,T].
\end{equation}
 if
 $\rho(\alpha^*,\alpha,T)<dist(y(0),bd\,{\mathcal{S}}).$
 \end{lemma}
 \begin{lemma}
\label{dop2}
 For a sequence of $\alpha_n\in{\ct{A}}$ and a
  sequence of  $\td{y}_n\in{\mm{Y}}[\alpha_n]$, let
 \[\rho(\alpha^*,\alpha_n,T)\to 0, y_n(0)\to \xi \textrm{ as } n\to\infty \]
 for some $T>0,$ $\xi\in int\,{\mathcal{S}}$.

 Then, the solutions
 $\td{y}_n$ converge to the solution of \rref{a} generated by $\alpha^*$ from the point $\xi$, and this convergence is uniform in $[0,T]$.
 \end{lemma}

 Hereinafter set $\Upsilon\rav U\times [1/2,\infty),$
  $\ct{A}=\ct{U}\times B(\mathbb{T},[1/2,\infty)),$ $\alpha^*=(u^*,1).$

  We will require the following property, which was essentially proved by Michel in \cite[Lemma]{Michel}:
\begin{lemma}
\label{mlemm}
 Consider Borel-measurable mappings $u\in\ct{U},v\in B(\mathbb{T},[1/2,\infty))$
 and the solutions of the system generated by them
\begin{eqnarray*}
    \dot{y}&=& v(t)\,f\big(y(t),u(t)\big),\qquad x(0)=b;\\
       \dot{z}&=& v(t),\qquad z(0)=0.
\end{eqnarray*}
 Then, there exists a control $u'\in\ct{U}$ and the trajectory $x'=x(b,u';\cdot)$ generated by it
  such that $x'(z(t))=x(t)$ for all $t\in[0,\tau_n]$ and
  \begin{equation}
  \label{532}
    \int_{0}^{\tau_n}v(t)e^{-rz(t)}f_0(x(t),u(t))\,dt=
    \int_{0}^{z(\tau_n)}e^{-rt}f_0(x'(t),u'(t))\,dt.
   \end{equation}
\end{lemma}
 {\doc}
  Note that every such map
 $z:\mm{T}\to\mm{T}$ is continuous,  strictly increasing, and reversible; denote the inverse map of
 $z$ by  $\zeta.$
 It would then suffice to set $u'(s)\rav u(\zeta(s)),x'(s)\rav x(\zeta(s))$ for all $s\leq z(\tau_n).$
 As proved in \cite[Lemm]{Michel}, in these circumstances, $y'=x(b_*,u';\cdot)$ and \rref{532}.
\bo

  For an unbounded sequence of positive numbers $\tau_n,$
   define the scalar function $h_n$ by
  the following rule:
$$ h_n(s)
  \rav e^{-rs}(J_{**}-J^0(b_*,0;\tau_n))\quad\forall s\in\mm{R}_{\geq 0}.$$
 Note that \rref{raz_} now implies
\begin{eqnarray}
  \label{razz}
  \displaystyle \lim_{n\to\infty} \sup_{s\in[-1,1]}{h}_n(s)=0
\end{eqnarray}

\begin{lemma}
\label{olemm}
   Suppose that $u^*$ is a uniformly overtaking optimal control of original problem \rref{sys0_}--\rref{sysK_}, i.e.
  in problem
  \begin{eqnarray*}
    \textrm{Minimize } l(b)+\int_{0}^\infty e^{-rt} f_0 (x,u)\, dt\\
    \textrm{subject to } \dot{x}=f(x,u),\quad u\in U,\\
    x(0)\in\ct{C}.
\end{eqnarray*}
  Assume the number $J_{**}$ to be validly defined by \rref{raz_}. Take an arbitrary unboundedly increasing sequence of times ~$\tau_n.$

  Assume that,
   for some unbounded sequence of positive numbers $\tau_n$,  a sequence of
functions $h_n\in C(\mm{R},\mm{R})$ satisfies \rref{razz}.

  Then,
   the sequence of optimal values of the problems
\begin{subequations}
\begin{eqnarray}
    {h}_n(z(\tau_n)-\tau_n)+l(x(0))+
    \int_{0}^{\tau_n}v(t)e^{-rz(t)}f_0(x(t),u(t))\,dt&\ & 
    \label{590__}\\
    \textrm{subject to }    \dot{x}= v(t)\,f\big(x(t),u(t)\big),\quad\dot{z}= v(t);&\ &\label{591__}\\
 \quad t>0,\quad u(t)\in U,\quad  |v(t)-1|\leq  e^{-t};&\ & \label{592__}
 \\
    x(0)\in\ct{C},\quad z(0)=0&\ & \label{593__}
\end{eqnarray}
 converges to $l(b_*)+J_{**}.$
\end{subequations}
\end{lemma}
{\doc }
 Note that the control $(u^*,1)$ is admissible for problem \rref{590__}--\rref{593__}, and, by the definition of $J_{**}$ and \rref{razz}, it provides the value of payoff that is arbitrary close to $l(b_*)+J_{**}$ (for large $n$).

By condition, there exists a sequence of positive $\omega_n$ that converges to zero such that ${h}_n(t)\leq \omega_n$  if $|t|<1$ for every $n\in\mm{N}$.

Assume the implication of the lemma to be false. Then, there exist a positive number $\epsi$, a sequence of initial conditions $b_n\in\ct{C}$,  and a sequence of controls $(u_n,v_n)$ with \rref{592_} such that, for any natural $n$,
the  trajectory  $(x_n,z_n)$ generated by the control $(u_n,v_n)$ from the position $(b_n,0)$ satisfies
\begin{eqnarray*}
l(b_n)&+&
 {h}_n(z_n(\tau_n)-\tau_n)+
 \int_{0}^{\tau_n}v_n(t)e^{-rz_n(t)}f_0(x_n(t),u_n(t))\,dt\leq
l(b_*)+J^{**}-4\epsi.
\end{eqnarray*}
Since we also have $|\dot{z}_n(\tau_n)-1|\leq e^{-t},$ we now know that
$|z_n(\tau_n)-\tau_n|<1,$ i.e. $|{h}_n(z(\tau_n)-\tau_n)|\leq \omega_n.$
Now, for all $n$ starting with a certain one,
$$l(b_n)+\int_{0}^{\tau_n}v_n(t)e^{-rz_n(t)}f_0(x_n(t),u_n(t))\,dt
\leq l(b_*)+J^{**}-3\epsi.$$

 Thanks to Lemma~\ref{mlemm}, for sufficiently large $n,$ there exists the control $u'_n\in\ct{U}$ and the trajectory  $x'_n=x(b_n,u'_n;\cdot)$ generated by it
 such that \rref{532} holds, whence
 $$l(b_n)+\int_0^{z_n(\tau_n)}e^{-rt}f_0(x'_n(t),u'_n(t))\,dt\leq  l(b_*)+J^{**}-2\epsi.$$
 Now, $z_n(\tau_n)\to\infty$ and \rref{raz_} imply that, for sufficiently large $n\in\mm{N}$,
  $$l(x'_n(0))+\int_0^{z_n(\tau_n)}e^{-rt}f_0(x'_n(t),u'_n(t))\,dt\leq l(b_*)+\int_0^{z_n(\tau_n)}e^{-rt}f_0(x^*(t),u^*(t))\,dt-\epsi.$$
 However, it contradicts the fact that $(x^*,u^*)$ is a uniformly overtaking optimal process for problem \rref{sys0_}--\rref{sysK_}.
 \bo

 This allows us to proceed to the actual proof of the main result.

\section{Proof of Theorem~\ref{4}.
 }
\label{doc}

\subsection{
{  { Choosing the metric}~$\rho$.}
}
  Consider the following system:
 \begin{subequations}
 \begin{eqnarray}
   \label{sys_x_n}
       \dot{x}&=& v\,f\big(x,u\big);\\
       \dot{z}&=& v;\label{sys_z_n}\\
       \dot{\psi}&=& -v\,\frac{\partial f}{\partial x}\big(x,u\big)
                    +\lambda v\,e^{-rz}\frac{\partial f_0}{\partial x}\big(x,u\big);\label{sys_psi_}\\
       \dot{\phi}&=& -\lambda r ve^{-rz}f_0(x,u);\label{sys_phi_n}\\
       \dot{\lambda}&=& 0. \label{sys_l_n}
\end{eqnarray}
\end{subequations}

  Remember that
   $\Upsilon=U\times[1/2,\infty).$
    Let $\Omega$ be a ball in ${\mathbb{X}}$ centered at $b_*$ with the radius $1/2.$
  Set
   $E\rav{\mathbb{X}}\times {\mathbb{R}}\times{\mathbb{X}}\times {\mathbb{R}}\times {\mathbb{R}},$ $y_*\rav(b_*,0,0,0,0)\in E.$
    Let ${\mathcal{S}}$ be a ball in $E$ centered at $y_*$ with the radius $2.$


Let the mapping $a:E\times \Upsilon\times{\mm{T}}\to E$
 be the right-hand side of system \rref{sys_x_n}--\rref{sys_l_n}.
  This system satisfies all the requirements we demand from a system \rref{a}.
  Consider mappings $w,\rho$ for such system~$a$
  with  designated control $\alpha^*=(u^*,1)$ and the compact set ${\mathcal{S}}$.

  Remember that ${\ct{A}}=\ct{U}\times B([1/2,\infty))$, and, for each $n\in{\mathbb{N}}$,
\[{\ct{A}}_{\tau_n}\rav \{\alpha=(u,v)\in{\ct{A}}\,|\,u(t)=u^*(t),v(t)=1 \ \forall t>\tau_n\}.\]
 By Lemma~\ref{dop0},   ${\ct{A}}_{\tau_n}$ is metrizable by  $(\alpha',\alpha'')\mapsto
  \rho(\alpha',\alpha'',\tau_n).$

\subsection{
  {Constructing the auxiliary optimal solution sequence}.}

 By Lemma~\ref{olemm}, there exists a sequence of positive numbers $\gamma_n$ converging to zero such that,  for any natural $n$, the optimal value for \rref{590__}--\rref{593__}
 is bounded from below by the value  $l(b_*)+J_{**}-\gamma_n^2.$
   Then, it is also a bound from below for the value of the following auxiliary minimum problem:
\begin{subequations}
\begin{eqnarray}
    l(x(0))+\int_{0}^{\tau_n}v(t)e^{-rz(t)}f_0(y(t),u(t))\,dt&+&h_n(z(\tau_n)-\tau_n)\nonumber\\
    +\gamma_n \rho(\alpha^*,\alpha,\tau_n)&+&\gamma_n||x(0)-b_*||
    \quad\label{590_}\\
    \textrm{subject to }    \dot{x}= v(t)\,f\big(x(t),u(t)\big),\quad\dot{z}= v(t);\label{591_}\\
 t\geq 0,\quad \alpha(t)=(u(t),v(t)),\quad u(t)\in U,\ |v(t)-1|\leq e^{-t}; \label{592_}
 \\
    x(0)\in\ct{C}\cap\Omega,\quad z(0)=0. \label{593_}
%
%
\end{eqnarray}
\end{subequations}
 Consider the set of all admissible controls $\alpha=(u,v)$ in this problem. This set contains $\alpha^*=(u^*,1),$ and is a  subspace
 of the complete metric space  ${\ct{A}}_{\tau_n}.$

    By the Ekeland principle \cite[Theorem~5.3.1]{ekeland},
  \cite[Theorem~2.1.3]{conv_new},
  for problem \rref{590_}--\rref{593_},
  there exists an  optimal pair $(b_n,\alpha_n)$  in the complete metric subspace of  $(\ct{C}\cap\Omega)\times{\ct{A}}_{\tau_n}$;
  denote by $(\td{x}_n,\td{z}_n)$ its  solution of \rref{591_}, \rref{593_}.
  Moreover (see \cite[Theorem~5.3.1,(i)]{ekeland},\cite[Theorem~2.1.3,(ii)]{conv_new}),
  \begin{subequations}
\begin{eqnarray}\nonumber
l(b_*)+{J}^0(b_*,0,u^*;\tau_n)+h_n(0)&\geq&\int_{0}^{\tau_n}v_n(t)e^{-rz_n(t)}f_0(\td{x}_n(t),\td{z}_n(t))\,dt\\
& & +l(\td{x}_n(0))+
h_n(\td{z}(\tau_n)-\tau_n)\nonumber\\
& &+\gamma_n \rho(\alpha_*,\alpha_n,\tau_n)+\gamma_n ||\td{x}_n(0)-b_*||,
\label{to_w_}\\
  ||\td{x}_n(0)-b_*||+\rho(\alpha_*,\alpha_n,\tau_n)&<&\gamma_n\to 0\textrm{ as }n\to\infty\label{1050}.
   \end{eqnarray}
   From \rref{591_} and \rref{592_} one can readily obtain $|\td{z}_n(\tau_n)-\tau_n|<1;$
   now, \rref{razz} implies
   \begin{eqnarray}
   \frac{dh_n}{ds}(\td{z}_n(\tau_n)-\tau_n)=-rh_n(\td{z}_n(\tau_n)-\tau_n)\to 0.\label{918}
\end{eqnarray}
    Let us show that
   \begin{eqnarray}
   \int_{0}^{\tau_n}v_n(t)e^{-rz_n(t)}f_0(\td{x}_n(t),\td{z}_n(t))\,dt\to J_{**}.\label{919}
\end{eqnarray}
\end{subequations}
    Indeed, to prove that the upper limit does not exceed $J_{**}$, it is sufficient to pass to the limit in \rref{to_w_} using  \rref{razz},\rref{1050}.
    On the other hand, as it was noted before, the integral can be estimated from below by the value $l(b_*)-l(b_n)+J_{**}-\gamma_n^2$ by virtue of Lemma~\ref{olemm}. However, the limit of this expression is also equal to  $J_{**}.$
    Thus, \rref{919} is proved.

Note that, by \rref{1050} and Lemma~\ref{dop0}, $\alpha_n=(u_n,v_n)$ converges in measure to $\alpha^*=(u^*,1)$ on the whole ${\mathbb{T}}.$
Passing to the subsequence if necessary, we can say that
$(u_n,v_n)$ converges to $(u^*,1)$ a.e. on ${\mathbb{T}}.$

\subsection{
  {Pontryagin Maximum Principle for auxiliary problem.}}

  Since $\alpha_n=(u_n,v_n)$ provides a minimum of problem \rref{590_}-\rref{593_}, it can, if need be, yield the Pontryagin Maximum Principle~\cite[Theorem~5.1.1]{cl_new}.
 Without loss of generality, we may \rref{1050} assume  that $\td{x}_n(0)\in int\,\Omega$ for all $n\in\mm{N}$. Then,
$N_L^{{\mathcal{C}}}(\td{x}_n(0))=N_L^{\Omega\cap{\mathcal{C}}}(\td{x}_n(0)).$

  Let    the function ${H}_n:{\mathbb{X}}\times{\mathbb{R}}\times{\Upsilon}\times{\mathbb{X}}\times{\mathbb{R}}\times{\mathbb{T}}\times {\mathbb{T}}\mapsto{\mathbb{R}}$
 be given by
  \[\displaystyle {H}_n(x,z,u,v,\psi,\phi,\lambda,t)\rav \psi v f(x,u)+\phi v  -\lambda v e^{-rz}f_0(x,u)-
  \lambda\gamma_n w(\alpha^*(t),(u,v),t).\]
  Then, by the  Maximum
Principle,
 there exist
 ${\lambda}_n\in(0,1]$, $\td{\psi}_n\in C({\mathbb{T}},{\mathbb{X}}),$
 $\td{\phi}_n\in C({\mathbb{T}},{\mathbb{R}})$
 with
 \begin{subequations}
\begin{eqnarray}\label{dob_}
\lambda+|\td{\phi}_n(0)|+||\td{\psi}_n(0)||=1
\end{eqnarray}
such that,
 for some $\zeta\in{\mathbb{X}} (||\zeta||\leq 1),$
 the transversality conditions
\begin{eqnarray}
\td\psi_n(0)&\in& \lambda_n\upartial_L^1 l(\td{x}_n(0))+\lambda_n\gamma_n\zeta+N_L^{{\mathcal{C}}}(\td{x}_n(0)),
\label{trans_0_max_}\\
-\td\phi_n(\tau_n)&=&\lambda_n\frac{dh_n}{ds}(\td{z}_n(\tau_n)-\tau_n),\label{trans_np_max_}\\
-\td\psi_n(\tau_n)&=&0\label{trans_n_max_}
   \end{eqnarray}
 hold,  and
 \begin{eqnarray}
\nonumber
-\dot{\td{\psi}}_n(t)&=&\frac{\partial
{H}_n}{\partial
x}\big(\td{x}_n(t),u_n(t),v_n(t),\td{\psi}_n(t),\td{\phi}_n(t),{\lambda}_n,t\big)
\\ \label{1575}
&=&v_n(t)\,\frac{\partial f}{\partial x}\big(\td{x}_n(t),u_n(t)\big)
                    \!\!-\!\!\lambda_n v_n(t)\,e^{-r\td{z}_n(t)}\frac{\partial f_0}{\partial x}\big(\td{x}_n(t),u_n(t)\big);\\
-\dot{\td{\phi}}_n(t)&=&\frac{\partial{H}_n}{\partial
z}\big(\td{x}_n(t),u_n(t),v_n(t),\td{\psi}_n(t),\td{\phi}_n(t),{\lambda}_n,t\big)\nonumber\\
&=&\lambda_n r v_n(t)e^{-r\td{z}_n(t)}f_0(\td{x}_n(t),u_n(t)\big);
\label{1675}\\
\nonumber
   \sup_{u'\in
        U,|v'-1|\leq e^{-t}}\!\!\!\!\!\!&\!&{H}_n\big(\td{x}_n(t),u',v',\td{\psi}_n(t),\td{\phi}_n(t),{\lambda}_n,t\big)\\
        &=&
           {H}_n\big(\td{x}_n(t),u_n(t),v_n(t),\td{\psi}_n(t),\td{\phi}_n(t),{\lambda}_n,t\big)
              \label{sys_max_}
 \end{eqnarray}
\end{subequations}
   also hold for a.e. $t\in[0,\tau_n]$.

\subsection{
  {Pontryagin Maximum Principle for overtaking optimal process}}
   Set $\td{y}_n\equiv (\td{y}_n,\td{z}_n,\td{\psi}_n,\td{\phi}_n,{\lambda}_n)$ for each $n\in{\mathbb{N}};$ note that this is a solution of \rref{sys_x_n}--\rref{sys_l_n}.

  By \rref{dob_},
  passing, if need be, to a subsequence, we can consider the subsequence of
  ${\lambda}_n \in (0,1]$ to tend to  some $\lambda^*\in[0,1]$ and a subsequence of $\big(\td{\psi}_n(0),\td{\phi}_n(0)\big)$ to converge to a certain $({\psi}^*_0,{\phi}^*_0)\in{\mathbb{X}}\times\mm{R}$ as well.
  We now have, for sufficiently large~$n$,
\begin{equation}\label{1208}
 ||\td{y}_n(0)-y_*||\leq \lambda^*+||\psi^*_0||+|\phi^*_0|+||\td{y}_n(0)-b_*||
 \leqref{dob_} 1+\gamma_n<2.
\end{equation}
   Thus,
  $\td{y}_n(0)\to(b_*,0,{\psi}^*_0,{\phi}^*_0,\lambda^*)\in int\,{\mathcal{S}}.$

   In addition, for each $T>0$, we have $\rho(\alpha^*,\alpha,T)\leq\rho(\alpha^*,\alpha,\tau_n)$ for all $u\in{\mathfrak{U}}$ if $T<\tau_n$; now, from  \rref{1050}, we have $\rho(\alpha^*,\alpha_n,T)\to 0.$ Therefore, by Lemma~\ref{dop2}, in every compact interval, the subsequence of $\td{x}_n$ uniformly converges to a solution $y^*$
  of system \rref{sys_x_n}--\rref{sys_l_n} generated by the control $(u^*,1)$, i.e.,
  the subsequence of $(\td{x}_n,\td{\psi}_n,\td{\phi}_n,\lambda_n)$ converges to the solution of
  \rref{sys_x_k}--\rref{sys_l_k}. Moreover, $y^*(0)=(b_*,0,\td{\psi}^*_0,\td{\phi}^*_0,\lambda^*).$
  Then, $y^*$ has the form
   $y^*(\cdot)=({x}^*(\cdot),\cdot,{\psi}^*(\cdot),{\phi}^*(\cdot),{\lambda}^*),$
   where functions ${\psi}^*,{\phi}^*$ are solutions of \rref{sys_psi_k} and of $\dot{\phi}^*=-\lambda^* r e^{-rt} f_0(x^*(t),u^*(t))$
   with initial conditions
   ${\psi}^*(0)={\psi}^*_0,$ ${\phi}^*(0)={\phi}^*_0.$

Remember that
$(u_n,v_n)$ converges a.a. to $(u^*,1)$. Then, $w((u^*(t),1),(u_n(t),v_n(t)),t)\to w((u^*(t),1),(u^*(t),1),t)=0$ for a.e. $t\in\mm{T}.$
 Now, passing to the limit in \rref{sys_max_}, we have,  for a.e. $t\in\mm{T}$,
 \begin{eqnarray}
 \sup_{u\in
        U,|v-1|\leq e^{-t}}
        \Big[\psi^*(t)vf\big({x}^*(t),u,t\big)+v\phi^*(t)-\lambda^* ve^{-rt}f_0\big({x}^*(t),u,t\big)\Big]=\nonumber\\
        \psi^*(t)f\big({x}^*(t),u^*(t),t\big)+\phi^*(t)-\lambda^* e^{-rt}f_0\big({x}^*(t),u^*(t),t\big).\ \label{337}
 \end{eqnarray}
 Setting $v=1$, we obtain  \rref{maxH} for  $({x}^*,{\psi}^*,{\lambda}^*)$ for almost every $t>0$.
 Thus, the limit $({x}^*,{\psi}^*,{\lambda}^*)$ satisfies system \rref{sys_x}--\rref{dob} for $u=u^*$, i.e., system \rref{sys_x_k}--\rref{sys_l_k}.


\subsection{
  {Backtracking}}

Since \rref{1050} and \rref{1208} imply 
$\rho(\alpha^*,\alpha_n,\tau_n)<\gamma_n<1/2<dist(\td{y}_n(0),bd\,{\mathcal{S}}),$
 and
$\td{y}_n\to y^*$, $\gamma_n\to 0$ as $n\to\infty$, we know that Lemma~\ref{dop} guarantees
\begin{equation}\label{1227}
   \varkappa(\td{y}_n(\tau_n),\tau_n)\to y^*(0).
 \end{equation}

From the position $\td{y}_n(\tau_n)$, launch in reverse time a solution $y_n$ of system
\rref{sys_x_n}--\rref{sys_l_n} with the help of the control $(u^*,1)$.
 Then, $y_n(0)=\varkappa(\td{y}_n(\tau_n),\tau_n)$ (see \rref{1667}).
   Note that $y_n=(x_n,z_n,\psi_n,\phi_n,\lambda_n)$ satisfies
   \rref{sys_x_k}--\rref{sys_l_k}, and
   $\psi_n(\tau_n)=\td{\psi}_n(\tau_n)=0$,
   $\phi_n(\tau_n)=\td{\phi}_n(\tau_n)=-\frac{dh_n}{ds}(\td{z}_n(\tau_n)-\tau_n)$.
    By the theorem on continuous dependence of the solution of a differential equation, \rref{1227} implies that
    the solution $y^*(\cdot)=\big(x^*(\cdot),\cdot,\lambda^*,\psi^*(\cdot),\phi^*(\cdot)\big)$  is the limit (in the compact-open topology) of
    $\td{y}_n.$

Note that the mappings $b\mapsto\upartial^{1}_L l(b)$, $b\mapsto N_L^{{\mathcal{C}}}(b)$ are upper semicontinuous; passing to the limit in
\rref{trans_0_max_}, we see that $\td\psi_n(0)\to\psi^*(0)$ and $\td{x}_n(0)\to x^*(0)=b_*$ imply \rref{400}.

 Since the supremum in \rref{337} contains a function that is linear in $v$ and that attains its maximum in $v$ at the interior point $v=1$, we have
$\psi^*(t)f\big({x}^*(t),u^*(t),t\big)+\phi^*(t)-\lambda^* e^{-rt}f_0\big({x}^*(t),u^*(t),t\big)=0$
 for almost every $t\in\mathbb{T}$,
i.e.,
$\phi^*(t)=-H(x^*(t),u^*(t),\psi^*(t),\lambda^*)$ for almost all  $t\in\mathbb{T}.$

   Now, $-H(x^*(t),u^*(t),\psi^*(t),\lambda^*)$ coincides with the limit of $\td{\phi}_n.$
Thanks to \rref{sys_phi_n} and \rref{trans_np_max_}, every $\td{\phi}_n$ also satisfies
 $$\dot{\td{\phi}}_n(t)=-\lambda_n r  e^{-r\td{z}_n(t)} f_0\big(\td{x}_n(t),u^*(t)\big),\quad
 -\td{\phi}_n(\tau_n)=\lambda_n\frac{dh_n}{ds}(\td{z}_n(\tau_n)-\tau_n).$$

 Then, for all $T\geq 0$,  the Hamiltonian $H^*[T]$ coincides with
 \begin{eqnarray*}
 \nonumber
 \lim_{n\to\infty}\Big[-\int_{T}^{\tau_n}\lambda_n r  e^{-r\td{z}_n(t)} f_0\big(\td{x}_n(t),u^*(t)\big)dt+\lambda_n\frac{dh_n}{ds}(\td{z}_n(\tau_n)-\tau_n)\Big]&\ravref{918}&\\
 \lim_{n\to\infty}\lambda_n r
 \bigg[-\int_{0}^{\tau_n} e^{-r\td{z}_n(t)}f_0\big(\td{x}_n(t),u^*(t)\big)\,dt
 +\int_{0}^{T} e^{-r\td{z}_n(t)} f_0\big(\td{x}_n(t),u^*(t)\big)\,dt\bigg]&\ravref{919}&\nonumber\\
  -\lambda^*r J_{**}+
 \lim_{n\to\infty}\lambda_nr \int_{0}^{T} e^{-r\td{z}_n(t)}f_0\big(\td{x}_n(t),u^*(t)\big)\,dt.\ \nonumber
\end{eqnarray*}
Passing to the limit, in view of the theorem on continuous dependence of the solution of a differential equation on its initial conditions,
we obtain
\begin{eqnarray}
 \nonumber H^*[T]
 &=&-\lambda^*rJ_{**}+\lambda^*r\int_{0}^{T}e^{-rt}f_0\big(x^*(t),u^*(t)\big)\,dt\\
 &=&\lambda^*r\big[-J_{**}+\bar{J}^0(b_*,u^*;T)\big].
 \label{304__}
\end{eqnarray}
Thus, \rref{412} is proved. Expression \rref{raz_} now implies \rref{413}.

  Note that although the constructed sequences converge to  $({\psi}^*_0,\phi^*,\lambda^*)$
  such that $||{\psi}^*_0||+|\phi^*_0|+{\lambda}^*=1$, we have
  $||{\psi}^*_0||+{\lambda}^*>0$. Indeed, we would otherwise have
  $\psi^*\equiv 0,$ $\lambda^*\equiv 0,$ whence $\phi^*(0)=1$ and
   $H(y^*(t),u^*(t),\psi^*(t),\lambda^*)\equiv 0,$ i.e. $H^*\equiv 0,$ which contradicts $H^*\equiv -\phi^*.$
   Thus, $||{\psi}^*_0||+{\lambda}^*>0$.
  Note that, since $||{\psi}_n(0)||\ravref{4B}\lambda_n||I(x_n(0),\tau_n)||,$
   ${\lambda}^*=0$ exactly when the sequence from  $||I(x_n(0),\tau_n)||$ is unboundedly increasing.

In case of    ${\lambda}^*>0, \lambda^*\neq 1$ note that
     relations  \rref{sys_psi},\rref{maxH},\rref{trans_0_max_}--\rref{1675} are preserved under multiplication of
    $({\psi}^*,\phi^*,\lambda^*)$ along with the subsequences $({\psi}_n,\phi_n,\lambda_n)$ by a positive number. Hence, by multiplying the triple
       $({\psi}^*,\phi^*,\lambda^*)$ along with the subsequences $({\psi}_n,\phi_n,\lambda_n)$ by the number
 $\frac{1}{{\lambda}^*},$
  we provide $\lambda^*=1$. Thus we can safely assume $\lambda^*\in\{0,1\}.$

  Expressions \rref{499},\rref{500} imply that, for all $T\geq 0$ for each  $n\in\mm{N}$,
    \begin{eqnarray*}
     \psi_n(T)&\in&\lambda_n\upartial_L (-\bar{J}^T)(x_n(T);u^*,\tau_n);\\
          (\psi_n(T), \lambda_n r J^T(x_n(T),0;\tau_n))&\in&
     \lambda_n\upartial^1_L (-J^T)({x}_n(T),0;\tau_n).
     \end{eqnarray*}
 Passing to the limit as $n\to\infty$, by virtue of $\lambda_n\to \lambda^*,$  ${x}_n\to x^*,$  ${\psi}_n\to \psi^*,$ and \rref{412}, we have \rref{409},\rref{410}.

  In case $\lambda^*=0$, \rref{304__} implies that $H^*\equiv 0,$ whence we obtain \rref{403};
  setting $\lambda^*=0$ in \rref{403}, we obtain \rref{404}.
\bo


\end{document}